\documentclass[11pt]{amsart} 
\usepackage{latexsym, amssymb}
\usepackage[nohug]{diagrams}
\diagramstyle[labelstyle=\scriptstyle]

\newtheorem{thm}{Theorem}[section]
\newtheorem{lemma}[thm]{Lemma}
\newtheorem{prop}[thm]{Proposition}
\newtheorem{cor}[thm]{Corollary}

\theoremstyle{definition}

\theoremstyle{remark}

\newcommand{\Z}{\mathbb{Z}}

\newcommand{\curly}[1]{\mathcal{#1}}

\newcommand{\B}{\curly{B}}

\def \<{\langle}
\def \>{\rangle}
\def \*Z {{{^*}\Z}}
\def \((  {(\!(}
\def \)) {)\!)}

\def \Ns{\operatorname{NSS}}

\numberwithin{equation}{section}

\allowdisplaybreaks[2]

\begin{document}

\title {Globalizing Locally Compact Local Groups} 

\author {Lou van den Dries and Isaac Goldbring}

\address {University of Illinois at Urbana-Champaign, Department of Mathematics,
1409 W. Green St., Urbana, IL 61801}
\email{vddries@illinois.edu}

\address {University of California, Los Angeles, Department of Mathematics, 520 Portola Plaza, Box 951555, Los Angeles, CA 90095-1555, USA}
\email{isaac@math.ucla.edu}
\urladdr{www.math.ucla.edu/~isaac}

\begin{abstract}
Every locally compact local group is locally isomorphic to a topological group.
\end{abstract}
\maketitle



\noindent
\section{Introduction}\label{in}

\noindent
In this paper a {\em local group\/} is a hausdorff topological space $G$ equipped with an element $1=1_G\in G$ (its identity) and continuous maps $x \mapsto x^{-1}: G \to G$ (inversion), and $(x,y) \mapsto xy: \Omega \to G$ (product) with
open $\Omega=\Omega_G\subseteq G\times G$, such that for all $x,y,z\in G$, \begin{enumerate}
\item $(1,x), (x,1)\in \Omega$ and $1x=x1=x$;
\item $(x, x^{-1}), (x^{-1},x)\in \Omega$ and $xx^{-1}=x^{-1}x=1$;
\item if $(x,y), (y,z)\in \Omega$ and either $(xy,z)\in \Omega$ or $(x,yz)\in \Omega$, then both  $(xy,z)$ and $(x,yz)$ belong to $\Omega$ and $(xy)z=x(yz)$.
\end{enumerate}
Also, ``topological group'' will stand for ``hausdorff topological group'', so
any topological group $G$ is a local group with $\Omega=G\times G$. From now on $G$ is a local group, and $a,b,c, x,y,z$ range over $G$. It follows easily that
if $xy=1$, then $x=y^{-1}$ and $y=x^{-1}$, and with a bit more effort, that
if $(x,y)\in \Omega$, then $(y^{-1}, x^{-1})\in \Omega$ and $(xy)^{-1}=y^{-1}x^{-1}$. Call $X\subseteq G$ {\em symmetric\/} if $X=X^{-1}$.
Given any symmetric open neighborhood $U$ of the identity of $G$ we have
a local group $G|U$; it has the subspace $U$ as underlying space, $1_G$ as its identity, the restriction of inversion to $U$ as its inversion, and the restriction of the product to 
$$\Omega_U:= \{(x,y)\in \Omega\cap (U\times U):\ xy\in U\}$$ as its product. 
Such a local group $G|U$ is called a {\em restriction\/} of $G$. 

A central question in the subject (see for example \cite{S1}) is the following: 
{\em When does a local group have a restriction that is also a restriction of a topological group?} 
We give here a complete answer in the locally compact case, as a consequence of the recent solution in \cite{Gold} of a local version of Hilbert's 5th problem and of earlier work by \`Swierczkowski~\cite{Sw}:

\begin{thm}\label{th1} If $G$ is locally compact, then some restriction of $G$ is a restriction of a topological group. 
\end{thm}

\noindent
Cartan \cite{Cart} established this for local Lie groups (local groups whose inversion and product are real analytic with respect to some real analytic manifold structure on the underlying space). This remains true for locally euclidean local groups, since these were shown in \cite{Gold} to be local Lie groups. 
The conclusion of the theorem is false for some local Banach-Lie groups; see \cite{vanEst}.

\medskip\noindent
The above notion of ``local group'' is that of \cite{Sw} and is
more strict than in \cite{Gold}, where inversion is only required to be defined on an open neighborhood of $1$, and the associativity axiom (3) is weaker. But the theorem goes through for the local groups of \cite{Gold}, since these have restrictions satisfying the local group axioms of the present paper; see \cite{Gold2}, Lemma 3.2.6. 

Note also that by local homogeneity \cite{Gold}, Lemma 2.16, the space 
$G$ is locally compact iff $1$ has a compact neighborhood. 
It is worth noting that in the totally disconnected case more is true (but the
result is much less deep):

\begin{prop}\label{th2} If $G$ is locally compact and totally disconnected, then some restriction of $G$ is a compact topological group.
\end{prop}

\noindent
This is proved just as in the global case for which we refer to \cite{MZ}, p. 54.  This proposition will not be used in the rest of the paper. 

\section{normal subgroups and strong morphisms}

\medskip\noindent
If $(a,b), (b,c)\in \Omega$ and $(ab,c)\in \Omega$, then we say that $abc$ is defined, and
set $$abc:= (ab)c=a(bc).$$ 
From now on $G'$ is also a local group. We take the local groups as the objects of a category;  a {\em morphism\/} $G \to G'$ in this category is a map
$\phi: G \to G'$ such that if $(x,y)\in \Omega_G$, then $(\phi(x), \phi(y))\in \Omega_{G'}$ and $\phi(xy)=\phi(x)\phi(y)$; composition of morphisms is given by
composition of maps. Below, morphisms are with respect to this category. We say that $G$ and $G'$ are {\em locally isomorphic\/} if they have restrictions that are isomorphic in this category. 

A {\em subgroup\/} of $G$ is a symmetric subset $H$ of $G$ such that 
$$1\in H, \quad H\times H\subseteq \Omega, \quad xy\in H \text{ for all }x,y\in H.$$
Note that such $H$ is a topological group in the obvious way. A {\em normal subgroup\/} of
$G$ is a subgroup $N$ of $G$ such that $axa^{-1}$ is defined and in $N$ for all $x\in N$ and $a\in G$ (in particular, $N\times G, G\times N\subseteq \Omega$). 

Suppose $N$ is a normal subgroup of $G$. It follows that for all $a,b$,
\begin{align*} aN=Na, \qquad &aN=bN\Leftrightarrow a\in bN,\\
\quad (a,b)\in \Omega\ \Longrightarrow\ &aN \times bN\subseteq \Omega  \text{ and }\ aN\cdot bN=(ab)N.
\end{align*} So the sets $aN$ form a partition of $G$. We make $G/N:= \{aN:\ a\in G\}$ into a local group (except possibly for the hausdorff requirement) by giving $G/N$ the quotient topology, taking $1N=N$ as the identity element, $aN \mapsto a^{-1}N$ as inversion, and setting
$$\Omega_{G/N}:= \{(aN, bN):\ (a,b)\in \Omega\}$$
with multiplication given by $(aN, bN) \mapsto (ab)N: \Omega_{G/N} \to G/N$. 
The canonical map $G \to G/N$ is open, and if $N$ is closed in $G$, then $G/N$ is hausdorff and this canonical map is a strong morphism in the following sense: 
a morphism $\phi: G \to G'$ is {\em strong\/} if
for all $(x,y)$ with $(\phi(x), \phi(y))\in \Omega_{G'}$ we have $(x,y)\in \Omega_G$.  If $\phi: G \to G'$ is a strong morphism, then its kernel
$\ker \phi:= \{x:\ \phi(x)=1_{G'}\}$ is a closed normal subgroup of $G$. 
The proof of the next lemma is obvious.

\begin{lemma}\label{L:strongrest}
If $\phi\colon G\to G'$ is a strong morphism and $V$ is a symmetric open neighborhood of the identity in $G'$, then $U:=\phi^{-1}(V)$ is a symmetric open neighborhood of the identity in $G$, and $\phi$ restricts to a strong morphism $G|U\to G'|V$.
\end{lemma}

\noindent
Sometimes a subgroup of $G$ is only ``partially'' normal, and we need to restrict
$G$ to get an actual normal subgroup. To describe this situation we say that a subgroup $N$ of $G$ has {\em normalizing neighborhood\/} $V$ if $V$ is a symmetric open neighborhood of $1$ in $G$ such that $N\subseteq V$ and for all $a\in V$ and $x\in N$, if $axa^{-1}$ is defined and in $V$, then $axa^{-1}\in N$. 

\begin{lemma}\label{c} Let $N$ be a compact subgroup of $G$ with normalizing neighborhood
$V$. Then there is an open symmetric neighborhood $U$ of $1$ in $G$ such that
$N\subseteq U \subseteq V$ and $N$ is a normal subgroup of $G|U$.
\end{lemma}
\begin{proof} By compactness of $N$ we can take a symmetric open neighborhood
$W$ of $1$ in $G$ such that for all $x\in N$ and
$a\in W$ we have $(x,a)\in \Omega$ and $xa\in V$, and $axa^{-1}$ is defined and in $V$. For $x\in N, a\in W$ we also have $(a,x)\in \Omega$ and
$ax\in V$, since $ax=(axa^{-1})a\in NW$. Thus for $U:= NW$ we have
$U=WN$. It follows easily that $U$ has the desired property.
\end{proof}

\begin{lemma}\label{d} Suppose $\iota : G \to G'$ is an injective open morphism. Then
$\iota$ restricts to an isomorphism $G|U \to G'|\iota U$, for some restriction $G|U$ of $G$.
\end{lemma}
\begin{proof} Any symmetric open neighborhood $U$ of $1$ in $G$ such that $U\times U\subseteq \Omega$ has the desired property.
\end{proof} 

\section{A theorem of \`Swierczkowski}

\noindent
The following theorem, slightly adapted to our situation, is from \textsection 11 of \cite{Sw}. It is closely related to earlier work by Smith~\cite{S2} and van Est~\cite{E}. 

\begin{thm}\label{Sw}
Let $L$ be a simply connected topological group whose second singular homology $H_2(L)$ vanishes. Let $V, \hat{V}, Q$ be symmetric open neighborhoods of the identity in $L$ such that
\begin{enumerate}
\item $\hat{V}\hat{V}\subseteq Q,\quad Q Q\subseteq V$,
\item $\hat{V}$ is connected, 
\item every closed curve in $\hat{V}\hat{V}$ is contractible in $Q$.
\end{enumerate}
If $\phi: G \to L|V$ is a surjective strong morphism, then there is an injective open morphism $G|\phi^{-1}\hat{V} \to G'$ into a topological group $G'$. \end{thm}

\noindent
In \cite{Sw} the topology on $G$ and $G'$ is discrete, which makes the continuity and openness of the morphisms trivial. Nevertheless, this discrete version yields the topological version above because we can assume
$G'$ to be generated by the image of $G$ so that the next well-known lemma is
applicable:

\begin{lemma}\label{L:top}
Suppose $H$ is a group and $\iota :G\to H$ is a morphism of the discrete local group underlying $G$ into the discrete group $H$, and suppose that $H$ is generated by 
$\iota(G)$.  Then there is a unique group topology on $H$ such that $\iota :G\to H$ is an open morphism.
\end{lemma}

\begin{proof}
Let $\B$ be the set of open neighborhoods of $1$ in $G$.  Let $\iota\B:=\{\iota(U) \ | \ U\in \B\}$.  It is easy to verify that $\iota\B$ is a neighborhood base at $1_{H}$ for a (necessarily unique) group topology on $H$, and that this topology has the desired properties.
\end{proof}

\noindent 
We need the following consequence of Theorem~\ref{Sw}:

\begin{cor}\label{T:enlarg} Suppose there exists a surjective strong morphism
$G \to L|V$ where $L$ is a Lie group and $V$ a symmetric open neighborhood of the identity in $L$. Then $G$ is locally isomorphic to a topological group.
\end{cor}

\begin{proof}
Shrinking $V$ and restricting $G$ accordingly we can replace $L$ by a locally isomorphic Lie group and so arrange that $L$ is simply connected. Then its second homotopy group $\pi_2(L)$ vanishes, \cite{Cart}, so by the Hurewicz Theorem of algebraic topology, \cite{Hatch}, Theorem 4-32, the singular homology group $H_2(L)$ also vanishes.
To obtain the desired conclusion, it suffices in view of Lemma~\ref{d} to find symmetric open neighborhoods $\hat{V}$ and $Q$ of the identity in $L$ satisfying (1)-(3) of Theorem~\ref{Sw}. Take a symmetric open neigborhood $Q$ of $1_L$ such that $QQ\subseteq V$. Lie groups are locally simply connected, so we can take a simply connected open neighborhood $Q'\subseteq Q$ of $1_L$. Next, take a symmetric open neighborhood $\hat{V}$ of $1_L$ such that $\hat{V}\hat{V}\subseteq Q'$. By passing to the connected component of $1_L$ in $\hat{V}$
we arrange that $\hat{V}$ is connected.  Every closed curve in $\hat{V} \hat{V}$ is a closed curve in $Q'$, so is contractible in $Q'$ and thus in $Q$.    
\end{proof}

\section{Proof of the Theorem}

\noindent
To apply \`Swierckowski's theorem to locally compact $G$, we need
a local analogue of a theorem of Yamabe. It appears in the thesis \cite{Gold2}, but not in \cite{Gold}, so we include here a proof using
\cite{Gold}. We say that $G$ has NSS (``no small subgroups'') if some neighborhood of $1$ contains no subgroup other than $\{1\}$. 

\begin{prop}\label{L:localyamabe}
Suppose $G$ is locally compact. Then $(G|U)/N$ has 
$\Ns$, for some compact normal subgroup $N$ of some restriction $G|U$ of $G$. 
\end{prop}

\begin{proof}
The inclusion diagram on the next page indicates the relevant symmetric subsets of $G$.  Lemma 9.3 of \cite{Gold} yields a compact subgroup $H$ of $G$ and a neighborhood $V$ of $1$ in $G$ such that every subgroup of $G$ contained in $V$ is contained in $H$. Since $H$ is compact, we have a compact normal subgroup $N$ of $H$ such that $N\subseteq H\cap V$ and $H/N$ has $\Ns$, by \cite{MZ}, page 99. This yields an open neighborhood $V'$ of $1$ in $G$ with $N\subseteq V'\subseteq V$ such that every subgroup of $G$ contained in $V'$ is 
contained in $N$.  Choose a symmetric open neighborhood $U$ of $1$ in $G$ so that $N\subseteq U\subseteq V'$ and all $abc$ with $a,b,c\in U$ are defined and
in $V'$. Then $U$ is a normalizing neighborhood of $N$ in $G$: if $a\in U$, then $aNa^{-1}$ is a subgroup of $G$ contained in $V'$, and thus $aNa^{-1}\subseteq N\subseteq U$.
By Lemma~\ref{c} we can shrink $U$ such that $N$ becomes a normal subgroup
of $G|U$. Then $(G|U)/N$ has $\Ns$.  
\end{proof}

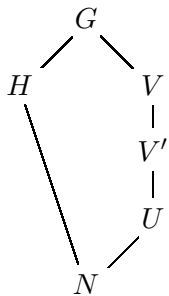
\begin{diagram}[size=1.15em]
 & &G & &\\
& \ldLine & &\rdLine &\\
H & & & &V\\ 
&\rdLine(2,6) & & &\dLine\\ 
& & & &V'\\
& & & &\dLine\\
& & & & U\\
& & & \ldLine&\\
& &N & &\\
\end{diagram}

\noindent
{\em Proof of Theorem~\ref{th1}}. Let $G$ be locally compact. Towards proving that
$G$ is locally isomorphic to a topological group we can replace $G$ by any restriction. After such a restriction, Proposition~\ref{L:localyamabe} gives a compact normal subgroup $N$ of $G$
such that $G/N$ has $\Ns$. Then by \cite{Gold}, \textsection 8, and Lemma~\ref{L:strongrest} we can arrange by a further restriction that $G/N$ is a restriction of a Lie group $L$.
It remains to apply Corollary~\ref{T:enlarg} to obtain that $G$ is locally isomorphic to a topological group.













\end{document}